\documentclass[preprint,12pt]{elsarticle}
\usepackage{graphicx,amsmath}

\newtheorem{defn}{Definition}

\begin{document}

\begin{frontmatter}
\title{Zero Dispersion and Zero Dissipation Implicit Runge-Kutta Methods for the Numerical Solution of Oscillating IVPs}

\author[UoP]{N. G. Tselios}
\ead{netselio@yahoo.gr}

\author[UoP]{Z. A. Anastassi}
\ead{zackanas@uop.gr}

\author[UoP]{T. E. Simos\fnref{simos}}
\fntext[simos]{Highly Cited Researcher, Active Member of the European Academy of Sciences and Arts, Address: Dr. T.E. Simos, 26 Menelaou Street, Amfithea - Paleon Faliron, GR-175 64 Athens, GREECE, Tel: 0030 210 94 20 091}
\ead{tsimos.conf@gmail.com, tsimos@mail.ariadne-t.gr}

\address[UoP]{Laboratory of Computer Sciences,\\
Department of Computer Science and Technology,\\
Faculty of Sciences and Technology, University of Peloponnese\\
GR-22 100 Tripolis, GREECE}

\begin{abstract}
In this paper we present two new methods based on an implicit
Runge-Kutta method Gauss which is of algebraic order fourth and has
two stages: the first one has zero dispersion and the second one has
zero dispersion and zero dissipation. The efficiency of these
methods is measured while integrating the radial Schr\"odinger
equation and other well known initial value problems.
\end{abstract}


\end{frontmatter}

\section{Introduction}

\quad We consider the radial Schr\"odinger equation:
\begin{equation}
\label{ode}
y''(x) = \left[\frac{µl(l + 1)}{ x^{2}} + V
(x)-E\right]y(x)
\end{equation}

 where $\frac{l(l+1)}{ x^{2}}$ is the
centrifugal potential, $V (x)$ is the potential, $E$ is the energy
and $W(x) = \frac{l(l+1)}{x^{2}} + V (x)$ is the effective
potential. It is valid that \[\lim_{\upsilon \rightarrow
\infty}{V(x)}=0\] and therefore \[\lim_{\upsilon \rightarrow
\infty}{W(x)}=0.\] We will study the case of $E
> 0$.

If we divide $[0,\infty]$ into small subintervals $[a_{i},b_{i}]$ so
that $W(x)$ is considered constant with value $\bar{W_{i}}$, then
the problem (1) is reduced to the approximation

\qquad $y''_{i} = \left(\bar{W}-E\right)y_{i}$, whose solution is

\begin{equation}
y_{i}(x) = A_{i}\,\ exp\left(\sqrt{\bar{W}-E} x\right)+ B_{i}\,\
exp\left(-\sqrt{\bar{W}-E} x\right),\quad A_{i}, B_{i} \in\Re.
\end{equation}

This form of Schr\"odinger equation shows why phase fittin is so
important when new methods are constructed. In the next section we
will present the most important parts of the theory used.

\quad The structure of the paper is as follows. Firstly in section
2, the basic theory of implicit Runge-Kutta methods is presented. In
section 3, the construction of the methods is introduced. In section
4 and 5, the calculation of the algebraic order and the
symplecticity respectively of the new methods is given. Then in
section 6, the numerical results of the new methods are presented
compared to classical RK methods from the literature while
integrating well known initial value problems and the radial
Schr\"odinger equation. Finally in section 7 our conclusions are
presented.

\section{Basic Theory}

\subsection{Implicit Method.}

The general form of an s-stage implicit Runge-Kutta method used for
the computation of the approximate value of $y_{n+1}(x)$ in Problem
\eqref{ode}, when $y_{n}(x)$ is known, is given from the following
procedure:

\begin{equation}
\label{rkmethod_gen0}
\begin{array}{l}
w_{i}=f(t_{n}+c_{i}h, y_{n}+h\sum\limits_{j=1}^{s}a_{ij}w_{j})\\
y_{n+1}=y_{n}+h\sum\limits_{i=1}^{s}b_{i}w_{i}
\end{array}
\end{equation}

when at least one $a_{ij} \neq 0$ exists with $i \leq j$.

An implicit Runge-Kutta method can also be presented using the
Butcher table below:

\begin{equation}
\label{rk_implicit}
\begin{array}{l|llllll}
  c_{1} & a_{11} & a_{12} &... & a_{1s-1} & a_{1s}\\
  c_{2} & a_{21} & a_{22} &... & a_{2s-1} & a_{2s}\\
  c_{3} & a_{31} & a_{32} &... & a_{3s-1} & a_{3s} \\
  \vdots&\vdots&\vdots&\ddots&\vdots&\vdots\\
  c_{s} & a_{s1} & a_{s2} &... & a_{s, s-1} & a_{ss}\\
  \hline
        & b_{1} & b_{2} &...  &b_{s-1} & b_{s}
\end{array}
\end{equation}

\subsection{Phase-Lag Analysis}

Let A and B, $s \times s$ matrices, be defined by $A=(a_{ij}), \quad
(1 \leq i, j \leq s$) and $B=(a_{ij}-b_{j}), \quad (1 \leq i, j \leq
s)$ respectively. When the method \eqref{rkmethod_gen0} is applied
to the linear equation

\begin{equation}
y'=qy, \quad q \in C
\end{equation}

the numerical solution is given by

\begin{equation}
y_{n+1}=P(z)y_{n}, \quad P(z)=\frac{det(I-zB)}{det(I-zA)}
\end{equation}

\noindent and can be written in the form

\begin{equation}
P(z)=K(\upsilon)+iL(\upsilon), \quad z=hq
\end{equation}

where $K(\upsilon)$ and $L(\upsilon)$ are functions of $\upsilon$
and $i=\sqrt{-1}$.

\begin{defn}
\label{defn1} \emph{ \cite{koto} In the implicit s-stage Runge-Kutta
method, presented in \eqref{rk_implicit}, the quantities
\begin{equation}
\phi(\upsilon)=\upsilon-arg(P(i\upsilon)), \quad
\alpha(\upsilon)=1-|P(i\upsilon)|, \quad \upsilon \in\ \Re
\end{equation}
are respectively called the phase-lag or dispersion error and the
dissipative error. If $\phi(\upsilon)=O(\upsilon^{q+1})$ and
$\alpha(\upsilon)=O(\upsilon^{r+1})$ then the method is said to be
of dispersive order $q$ and dissipative order $r$. }
\end{defn}

\subsection{Stability}

\begin{defn}
\label{defn2} \emph{ \cite{butcher} The stability function for an
implicit Runge-Kutta method is the rational function
\begin{equation}
 R(z)=\frac{det(I-zA+zeb^{T})}{det(I-zA)},
\end{equation}
where the vector $e=(1, ..., 1)^{T}$, and that a method is A-stable
if $|R(z)| \leq 1$, \quad whenever $Re(z) \leq 0$, where $Re(z)$ is
the real part of $z$. }
\end{defn}

\section{Construction of the new Runge-Kutta methods}

We consider the implicit Runge-Kutta method of Gauss, which is of
algebraic order fourth and has two stages. The coefficients are
shown in Table \ref{rkmethod_1}.

\begin{center}
\begin{equation}
\label{rkmethod_1}
\begin{array}{l|lll}
  \frac{1}{2}-\frac{\sqrt{3}}{6} & \frac{1}{4}& \frac{1}{4}-\frac{\sqrt{3}}{6} \\
  \frac{1}{2}+\frac{\sqrt{3}}{6}& \frac{1}{4}+\frac{\sqrt{3}}{6}\ \ \  & \frac{1}{4}\\

  \hline
      & \frac{1}{2} & \frac{1}{2}
\end{array}
\end{equation}
\end{center}
Below we present the construction of the methods.

\subsection{Construction of the new method with zero phase lag}

We consider all the values of Table \ref{rkmethod_1} except $b_2$.
By evaluating the phase-lag of this method, defined in Definition
\ref{defn1}, and by solving $\phi(v)=0$ towards $b_{2}$, the result
is:

\begin{equation}
\label{newmethod1}
\begin{array}{l}
\displaystyle b_{2}=\frac{1}{2}\,{\frac
{-6\,{\,\upsilon}^{3}+6\,{\,\upsilon}^{3}\,\sqrt {3}+72\,
\,\upsilon-\,\mathrm{tan} \left( \,\upsilon \right)
{\,\upsilon}^{4}+24\,\,\mathrm{tan} \left(\,\upsilon \right)
{\,\upsilon}^{2}}{ \,\mathrm{tan} \left( \,\upsilon \right) \,\sqrt
{3}{\,\upsilon}^{4}+6\,{\,\upsilon}^{3}+6\,{ \,\upsilon}^{3}\,\sqrt
{3} }} \\ \displaystyle \qquad \frac{+\tan \left( \upsilon \right)
{\upsilon }^{4}\sqrt {3}-12\,\mathrm{tan} \left( \upsilon \right)
{\upsilon}^{2}\sqrt {3 }-144\,\mathrm{tan} \left( \upsilon
\right)}{-36\,\upsilon^{2}\,\mathrm{tan} \left( \upsilon \right) -12
\,\upsilon^{2}\,\mathrm{tan} \left( \upsilon \right) \sqrt {3}-72
\,\upsilon }
\end{array}
\end{equation}

The Taylor series expansion of $b_{2}$ is shown below:

\begin{eqnarray*}
\begin{array}{l}
\displaystyle b_{2_{taylor}}= \ \frac{1}{2} \ + \ {\frac
{1}{720}}\,{\upsilon}^{4} \ + \  \left( {\ \frac {1}{6720}} \ - \
{\frac {1}{8640}}\,\sqrt {3} \  \right) {\upsilon}^{6} \ \ +\ldots
\end{array}
\end{eqnarray*}

In the last equation we observe that \[\lim_{\upsilon \rightarrow
0}{b_{2_{taylor}}}=\frac{1}{2}\] namely when the step-length tends
to zero the coefficient of the method Gauss appears.

\subsection{Construction of the new method with zero phase lag and zero
dissipation}

We consider all the values of Table \ref{rkmethod_1} except two:
$b_{2}$ and $a_{22}$. An extra equation (apart from the equation of
the phase lag) must hold, in order to achieve zero phase-lag and
zero dissipation. The two equations are $\phi(v)=0$ and
$\alpha(v)=0$.

After satisfying the above two equations, by solving towards $b_{2}$
and $a_{22}$, the result is:

\begin{equation}
\begin{array}{l}
\displaystyle {b_{2}}= \frac{1}{6}\cdot\frac{A}{B}
\end{array}
\end{equation}

where

\begin{eqnarray*}
\displaystyle A & = & -24\,\sin \left( \upsilon
\right) \upsilon\,\sqrt {3}+2 \,\sin \left( \upsilon \right)
{\upsilon}^{3}\sqrt {3}+12\,{\upsilon}^ {2}\cos \left( \upsilon
\right) \nonumber \\ &  & - 36\,\sin \left( \upsilon \right)
\upsilon-3\,\sin \left( \upsilon \right) {\upsilon}^{3}-144\,\cos
\left( \upsilon \right) \nonumber \\ &  & +12\,\sqrt {-{\frac
{-4\,{\upsilon}^{4}+2\,{ \upsilon}^{4}\sqrt
{3}+24\,{\upsilon}^{2}-24\,{\upsilon}^{2}\sqrt {3}- 144}{ \left(
\cos \left( \upsilon \right)  \right) ^{2}}}}\cos \left( \upsilon
\right)
\end{eqnarray*}

\begin{eqnarray*}
\displaystyle B={\upsilon} \left( 4\,\upsilon\,\sqrt {3 }\cos \left(
\upsilon \right) +\sin \left( \upsilon \right) {\upsilon}
^{2}-12\,\sin \left( \upsilon \right) \right)
\end{eqnarray*}

and

\begin{equation}
\begin{array}{l}
\displaystyle {a_{22}} = -\frac{1}{12}\cdot\frac{C}{D}
\end{array}
\end{equation}

with

\begin{eqnarray*}
\displaystyle C & = & 41472+72\,{\upsilon}^{4}\sqrt
{3}+864\,{\upsilon}^{2}- 90\,{\upsilon}^{6}\sqrt
{3}+1152\,{\upsilon}^{4}+  162\,{\upsilon}^{6} \nonumber  \\ & &  +
12096\,{\upsilon}^{2}\sqrt {3}-14688\,{\upsilon}^{2} \left( \cos
\left( \upsilon \right)  \right) ^{2} -144\,{\upsilon}^{4} \left(
\cos \left( \upsilon \right)  \right) ^{2} \nonumber \\ & & +41472\,
\left( \cos \left( \upsilon \right)  \right) ^{2}-120\,\sin \left(
\upsilon \right) { \upsilon}^{5}\sqrt {3}\cos \left( \upsilon
\right) \nonumber \\ & & -864\,\sin \left( \upsilon \right)
\upsilon\, T_{0}\cos \left( \upsilon \right) -504\, T_{0}\sin \left(
\upsilon \right) {\upsilon}^{3} }{ \left( \cos \left( \upsilon
\right) \right) ^{2}\cos \left( \upsilon \right) \nonumber \\ & & -{
\upsilon}^{7}\sin \left( \upsilon \right) \sqrt {3}\cos \left(
\upsilon \right) +20736\,\sin \left( \upsilon \right)
\upsilon\,\sqrt {3}\cos \left( \upsilon \right) \nonumber \\ & &
-1872\,\sin \left( \upsilon \right) { \upsilon}^{3}\sqrt {3}\cos
\left( \upsilon \right) -288\, T_{0} \left( \cos \left( \upsilon
\right)  \right) ^{2}{\upsilon}^{2}\sqrt {3} \left( \cos \left(
\upsilon \right) \right) ^{2} \nonumber \\ & &
+48\,{\upsilon}^{4}T_{0} \left( \cos \left( \upsilon \right) \right)
^{2}\sqrt {3} \left( \cos \left( \upsilon \right) \right)
^{2}+576\,{\upsilon}^{2}T_{0} \left( \cos \left( \upsilon \right)
\right) ^{2} \left( \cos \left( \upsilon \right) \right) ^{2}
\nonumber
\\ & & -36\,{\upsilon}^{5}\sin \left(
\upsilon \right) \cos \left( \upsilon \right) +10368\,\sin \left(
\upsilon \right) \upsilon\,\cos \left( \upsilon \right) \nonumber
\\& & +1728 \,\sin \left( \upsilon \right) {\upsilon}^{3}\cos \left(
\upsilon \right) +3\,{\upsilon}^{7}\sin \left( \upsilon \right) \cos
\left( \upsilon \right) \nonumber \\ & & -36\,{\upsilon}^{4}T_{0}
\left( \cos \left( \upsilon \right)  \right) ^{2}+132\,\sin \left(
\upsilon \right) {\upsilon}^{3}\sqrt {3}T_{0} \cos \left( \upsilon
\right) \nonumber \\& & +{\upsilon}^{5}\sin \left( \upsilon \right)
\sqrt {3}T_{0}\cos \left( \upsilon \right) -1728\,\sin \left(
\upsilon \right) \upsilon\,\sqrt {3}T_{0}\cos \left( \upsilon
\right) \nonumber \\& & +18\,{\upsilon}^ {6} \left( \cos \left(
\upsilon \right) \right) ^{2}-6912\,T_{0} \left( \cos \left(
\upsilon \right) \right) ^{2} 6\,{\upsilon}^{6}\sqrt {3} \left( \cos
\left( \upsilon \right) \right) ^{2} \nonumber \\& &
-1728\,{\upsilon}^{2}\sqrt {3} \left( \cos \left( \upsilon \right)
\right) ^{2}-360\,{\upsilon}^{4}\sqrt {3} \left( \cos \left(
\upsilon \right)  \right) ^{2}
\end{eqnarray*}

\begin{eqnarray*}
\displaystyle D & = & -144\,{\upsilon}^{4} \left( \cos \left(
\upsilon \right)  \right) ^{2}-6\,{\upsilon}^{6} \left( \cos \left(
\upsilon \right)  \right) ^{2}+576\,{\upsilon}{2}\,\sqrt {3}
+864\,{\upsilon}^{2} \nonumber \\& &-144\, \sin \left( \upsilon
\right) {\upsilon}^{3}\sqrt {3}\cos \left( \upsilon \right)
-2\,{\upsilon}^{6}\sqrt {3}-6\,{\upsilon}^{5}+120\,{ \upsilon}^{4}
\left( \cos \left( \upsilon \right)  \right) ^{3}\sqrt { 3}
\nonumber
\\& &-36\,{\upsilon}^{5}\cos \left( \upsilon \right) \sin \left(
\upsilon \right) -1728\,{\upsilon}^{2}\,\sqrt {3} \left( \cos \left(
\upsilon \right)  \right) ^{2} \nonumber \\& &-24\,\sin \left(
\upsilon \right) {\upsilon}^{5} \sqrt {3}\cos \left( \upsilon
\right) -1152\,{\upsilon}^{3}\cos \left( \upsilon \right) \sin
\left( \upsilon \right) \nonumber \\& & +3456\,\cos \left( \upsilon
\right) \sin \left( \upsilon \right)\,{\upsilon}
+96\,{\upsilon}^{2}\,\sqrt {3} \left( \cos \left( \upsilon \right)
\right) ^{2}-{\upsilon}^{7}\sin \left( \upsilon \right) \sqrt
{3}\cos \left( \upsilon \right) \nonumber
\\& &+24T_{0}\,\sin \left( \upsilon \right) {\upsilon}^{3}\cos \left( \upsilon
\right) +12\,{\upsilon}^{4} \left( \cos \left( \upsilon \right)
\right) ^{2} \nonumber \\& &+14\,{\upsilon}^{6} \sqrt {3} \left(
\cos \left( \upsilon \right)  \right) ^{2}+3\,{ \upsilon}^{7}\sin
\left( \upsilon \right) \cos \left( \upsilon \right) \nonumber \\& &
+{\upsilon}^{6}\sin \left( \upsilon \right) \sqrt {3}\cos \left(
\upsilon \right) -12\,\sin \left( \upsilon \right)
{\upsilon}^{4}\sqrt {3}\cos \left( \upsilon \right) \nonumber\\& &
-288\,\sin \left( \upsilon \right) \cos \left( \upsilon \right)
-24\,{\upsilon}^{3}\sqrt {3}-864\,\upsilon\, \left( \cos \left(
\upsilon \right)  \right) ^{2}
\end{eqnarray*}

where
\begin{eqnarray*}
\displaystyle T_{0}=\sqrt {-{\frac {2\,{\upsilon}^{4}\sqrt
{3}-4\,{\upsilon}^{4}+24\,{ \upsilon}^{2}-24\,{\upsilon}^{2}\sqrt
{3}-144}{ \left( \cos \left( \upsilon \right)  \right) ^{2}}}}
\end{eqnarray*}

The Taylor series expansion of $b_{2}$ and $a_{22}$ are shown below:

\begin{eqnarray*}
\displaystyle b_{2_{taylor}}={\frac {1}{2}}+{\frac
{1}{720}}\,{\upsilon}^{4}+{\frac {1}{10080}}\,{ \frac {5\,\sqrt
{3}-8}{-3+\sqrt {3}}}\,{\upsilon}^{6}+ \ldots
\end{eqnarray*}

\begin{eqnarray*}
\displaystyle a_{22_{taylor}}={\frac
{1}{4}}+{\frac {1}{2160}}\,{\frac {5\,\sqrt {3}-9}{-2+\sqrt {3
}}}\,{\upsilon}^{4}-{\frac {1}{181440}}\,{\frac {220\,\sqrt
{3}-381}{
 \left( -2+\sqrt {3} \right) ^{2}}}\,{\upsilon}^{6}+ \ldots
\end{eqnarray*}

In the last equations we observe that the limits when
$\upsilon\rightarrow 0$ are equal to the corresponding coefficients
of the Gauss method.

\section{Algebraic order of the new methods}
\label{sec4}

The following 8 equations must be satisfied so that the new method
maintains the fourth algebraic order of the corresponding classical
method presented in Table \ref{rkmethod_1}. The number of stages is
symbolized by $s$, where $s=2$.

\begin{equation}
\begin{array}{cc}\label{alg5_gen}
\textbf{1st Alg. Order (1 equation)}\\
{\sum\limits_{i = 1}^{s} {b_{i}} }  = 1
\end{array}
\end{equation}

\begin{equation}
\begin{array}{cc}
\textbf{2st Alg. Order (2 equations)}\\
{\sum\limits_{i = 1}^{s} {b_{i}} } c_{i} = \frac{1}{2}
\end{array}
\end{equation}

\begin{equation}
\begin{array}{cc}
\textbf{3st Alg. Order (4 equations)}\\
{\sum\limits_{i = 1}^{s} {b_{i}} } c_{i} ^{2} = \frac{1}{3}\\
{\sum\limits_{i,j = 1}^{s} {b_{i}} } a_{ij} c_{j} = \frac{1}{6}
\end{array}
\end{equation}

\begin{equation}
\begin{array}{cc}
\textbf{4st Alg.Order (8 equations)}\\
{\sum\limits_{i = 1}^{s} {b_{i}} } c_{i} ^{3} = \frac{1}{4}\\
{\sum\limits_{i,j = 1}^{s} {b_{i}} } c_{i} a_{ij} c_{j} = \frac{1}{8}\\
{\sum\limits_{i,j = 1}^{s} {b_{i}} } a_{ij} c_{j}^2 = \frac{1}{12}\\
{\sum\limits_{i,j,k = 1}^{s} {b_{i}} } a_{ij} a_{jk} c_{k} =
\frac{1}{24}
\end{array}
\end{equation}

\subsection{Remainders for the first method (algebraic conditions)}

We present the remainders of the eight equations, that is the
difference of the right part minus the left part, for the first
method:

\begin{equation}
\begin{array}{l}
rem_{1}={\frac {1}{720}}\,\, {\upsilon}^{4}+ \left( {\frac
{1}{6720}}-{\frac {1}{ 8640}}\,\sqrt {3} \right) {\upsilon}^{6}
+\ldots\\
rem_{2}=\left({\frac {1}{1440}}\,+{\frac {\sqrt {3}}{4320}}\,\,
\right){\upsilon}^{4} +\, \left({\frac {1}{60480}}\,-{\frac {\sqrt
{3}}{30240}}\,\, \right){ \upsilon}^{6} + \ldots\\
rem_{3}= \left({\frac {1}{2160}}\,+{\frac {\sqrt {3}}{4320}}\,\,
\right){\upsilon}^{4}-\, \left({\frac {1}{120960}}\,-{\frac
{\sqrt {3}}{72576}}\,\,\right){ \upsilon}^{6}+ \ldots \\
rem_{4}= \left({\frac {1}{4320}}\,+{\frac {\sqrt {3}}{8640}}\,\,
\right){\upsilon}^{4}-\, \left({\frac {1}{241920}}\,-{\frac {\sqrt
{3}}{145152}}\,\, \right){ \upsilon}^{6}+ \ldots \\
rem_{5}= \left({\frac {1}{2880}}\,+{\frac {\sqrt {3}}{5184}}\,\,
\right){\upsilon}^{4}-\, \left({\frac {1}{90720}}\,-{\frac {\sqrt
{3}}{120960}}\,\, \right){ \upsilon}^{6}
+ \ldots \\
rem_{6}=\, \left({\frac {1}{5760}}\,+{\frac {\sqrt{3}}{10368}}\,\,
\right){\upsilon}^{4}-\, \left({\frac {1}{181440}}\,-{\frac
{\sqrt{3}}{241920}}\right)\,{ \upsilon}^{6}
+ \ldots \\
rem_{7}=\, \left({\frac {1}{8640}}\,+{\frac {\sqrt
{3}}{12960}}\,\,\right){\upsilon}^{4}-\,\left({\frac
{1}{145152}}\,-{\frac {\sqrt {3}}{725760}}\right)\,\,{ \upsilon}^{6}
+ \ldots \\
rem_{8}= \left({\frac {1}{17280}}\,+{\frac {\sqrt{3}}{25920}}\,\,
\right){\upsilon}^{4}-\, \left({\frac {1}{290304}}\,-{\frac
{\sqrt{3}}{1451520}}\,\, \right){ \upsilon}^{6}
+ \ldots \\
\end{array}
\end{equation}

We see that the eight equations are held, when $h\rightarrow 0
\Rightarrow v\rightarrow 0$. This means that the new method
maintains the algebraic order of the corresponding classical method.

\subsection{Remainders for the second method (algebraic conditions)}

Now we present the remainders of the equations for the second
method:

\begin{equation}
\begin{array}{l}
rem_{1}={\frac {1}{720}}\, {\upsilon}^{4}+{\frac {1}{10080}}\,{\frac
{-8+5\, \sqrt {3}}{\sqrt {3}-3}}\, {\upsilon}^{6}+ \ldots \\
rem_{2}=\left( {\frac {1}{1440}}+{\frac {1}{4320}}\,\sqrt {3}
\right) { \upsilon}^{4}+{\frac {1}{60480}}\,{\frac { \left(
-8+5\,\sqrt {3} \right)  \left( 3+\sqrt {3} \right) }{\sqrt
{3}-3}}\,\,{\upsilon}^{6}
+ \ldots \\
rem_{3}={\frac {1}{25920}}\, \left( 3+\sqrt {3} \right)
^{2}{\upsilon}^{4}+{ \frac {1}{362880}}\,{\frac { \left( -8+5\,\sqrt
{3} \right)  \left( 3+ \sqrt {3} \right) ^{2}}{\sqrt
{3}-3}}\,\,{\upsilon}^{6}+ \ldots \\
rem_{4}=-{\frac {1}{8640}}\,{\frac {11\,\sqrt {3}-21}{ \left( \sqrt
{3}-3 \right)  \left( -2+\sqrt {3} \right)
}}\,\,{\upsilon}^{4}+{\frac {1}{ 181440}}\,{\frac {137\,\sqrt
{3}-237}{ \left( \sqrt {3}-3 \right)
\left( -2+\sqrt {3} \right) ^{2}}}\,\,{\upsilon}^{6}+ \ldots \\
rem_{5}={\frac {1}{155520}}\, \left( 3+\sqrt {3} \right)
^{3}{\upsilon}^{4}+{ \frac {1}{2177280}}\,{\frac { \left(
-8+5\,\sqrt {3} \right)  \left( 3 +\sqrt {3} \right) ^{3}}{\sqrt
{3}-3}}\,\,{\upsilon}^{6}+ \ldots \\
rem_{6}=-{\frac {1}{51840}}\,{\frac { \left( 3+\sqrt {3} \right)
\left( 11\, \sqrt {3}-21 \right) }{ \left( \sqrt {3}-3 \right)
\left( -2+\sqrt {3 } \right) }}\,\,{\upsilon}^{4}+{\frac
{1}{1088640}}\,{\frac { \left( 3+ \sqrt {3} \right)  \left(
137\,\sqrt {3}-237 \right) }{ \left( \sqrt { 3}-3 \right)  \left(
-2+\sqrt {3} \right) ^{2}}}\,\,{\upsilon}^{6}+ \ldots \\
rem_{7}=-{\frac {1}{8640}}\, \left( -2+\sqrt {3} \right)
^{-1}{\upsilon}^{4}+ {\frac {1}{725760}}\,{\frac {-53+31\,\sqrt
{3}}{ \left( \sqrt {3}-3 \right)  \left( -2+\sqrt {3} \right) ^{2}}}\,\,{\upsilon}^{6}+ \ldots \\
rem_{8}=-{\frac {1}{17280}}\,{\frac {3385\,\sqrt {3}-5863}{ \left(
\sqrt {3}- 3 \right)  \left( -2+\sqrt {3} \right)
^{5}}}\,\,{\upsilon}^{4}+{\frac {1} {290304}}\,{\frac {28121\,\sqrt
{3}-48707}{ \left( \sqrt {3}-3 \right)  \left( -2+\sqrt {3} \right) ^{6}}}\,\,{\upsilon}^{6}+ \ldots \\
\end{array}
\end{equation}
We see that for $\upsilon=0$ the eight equations are held for this
method too. Thus the new method has also fourth algebraic order.

\section{Symplecticity of the new methods}

\textbf{Theorem}\quad The Runge-Kutta method (3)-(4) is symplectic
when the following equalities are satisfied

\begin{equation}
b_{i}a_{ij}+b_{j}a_{ji}=b_{i}b_{j},  \qquad 1\leq i,j\leq s.
\end{equation}

As a classical example we mention the Gauss methods as symplectic
Runge-Kutta methods. It should be noted that symplectic Runge-Kutta
methods are always implicit.

Thus according to the above theorem, the three equations must be
satisfied so that the symplecticity of the new methods will be
maintained.

\begin{equation}
b_{1}a_{11}+b_{1}a_{11}=b_{1}b_{1}
\end{equation}

\begin{equation}
b_{2}a_{22}+b_{2}a_{22}=b_{2}b_{2}
\end{equation}

\begin{equation}
b_{1}a_{12}+b_{2}a_{21}=b_{1}b_{2}
\end{equation}

\subsection{Remainders for the first method (simplecticity conditions)}

We present the remainder of the three equations, that is the
difference of the right part minus the left part, for the first
method:

\begin{equation}
\begin{array}{l}
rem_{1}= 0\\
rem_{2}=-{\frac {1}{1440}}\,{\upsilon}^{4}+ \left( -{\frac
{1}{13440}}+{\frac { 1}{17280}}\,\sqrt {3} \right){\upsilon}^{6}+ \ldots \\
rem_{3}= \left( -{\frac {1}{2880}}+{\frac {1}{4320}}\,\sqrt {3}
\right){ \upsilon}^{4}+ \left( -{\frac {23}{241920}}+{\frac
{13}{241920}}\, \sqrt {3} \right){\upsilon}^{6}+ \ldots \\
\end{array}
\end{equation}

We see that for the three equations are held, when $h\rightarrow 0
\Rightarrow v\rightarrow 0$. That means that the new method
maintains the symplecticity  of the corresponding classical method.

\subsection{Remainders for the second method (simplecticity conditions)}

Now we present the remainders of the equations for the second
method:

\begin{equation}
\begin{array}{l}
rem_{1}= 0\\
rem_{2}={\frac {1}{720}}\,{\frac {-45+26\,\sqrt {3}}{ \left(
-2+\sqrt {3} \right)  \left( -3+\sqrt {3} \right)
^{2}}}\,{\upsilon}^{4}+{\frac {1}{ 10080}}\,{\frac {545\,\sqrt
{3}-944}{ \left( -2+\sqrt {3} \right) ^{2}\left( -3+\sqrt {3} \right) ^{3}}}\,{\upsilon}^{6}+ \ldots \\
rem_{3}=-{\frac {1}{2880}}\,{\frac {-5+3\,\sqrt {3}}{-3+\sqrt
{3}}}\,{\upsilon} ^{4}-{\frac {1}{120960}}\,{\frac {-54+31\,\sqrt
{3}}{-3+\sqrt {3}}}\,{ \upsilon}^{6}+ \ldots \\
\end{array}
\end{equation}

We see that for $\upsilon=0$ the three equations are held for this
method too. Thus the new method is also symplectic.

\section{Numerical Results}

\subsection{The methods}

In order to measure the efficiency of the methods constructed in
this paper we compare them to some already known methods, presenting
the results of the best six.

\textbf{I.} \quad \ Method \textbf{G2-PL-D} constructed in this
paper, where G2-PL-D means the method Gauss two-stages, fourth-order
with zero phase-lag and zero dissipation.

\textbf{II.} \quad Method \textbf{G2-PL} constructed in this paper,
where G2-PL means the method Gauss two-stage, fourth-order with zero
phase-lag.

\textbf{III.} \ \ Method \textbf{G2}: The classical two-stages and
fourth-order Gauss method (see \cite{hairer}).

\textbf{IV.}\quad Method \textbf{SDIRK(3,6,3)}: The Singly
Diagonally-Implicit Runge-Kutta method of J. M. Franco, I. Gomez, L.
Randez, is third-stage, third algebraic order,sixth dispersive order
and third dissipative order (see \cite{franco}).

\textbf{V.}\ \ \ \ \ Method \textbf{Radau I}: The classical third
order Radau method (see \cite{butcher}).

\textbf{VI.}\quad Method \textbf{Lobatto IIIC}: The classical
fourth order Lobatto method (see \cite{butcher}).

\subsection{The Problems}

\subsubsection{Inverse Resonance Problem}

The efficiency of the two new constructed methods will be measured
through the integration of problem \eqref{ode}  with $l = 0$ at the
interval $[0,15]$ using the well known Woods-Saxon potential

\begin{equation}
\begin{array}{l}
\displaystyle V(x)=\frac{u_{0}}{1+q} + \frac{u_{1}q}{(1+q)^{2}},
\quad q=exp \left( \frac{x-x_{0}}{a} \right),
\end{array}
\end{equation}

where \quad $u_{0}=-50,\quad a=0.6,\quad x_{0}=7\quad$ and $\quad
\displaystyle u_{1}=-\frac{u_{0}}{a}$

and with boundary condition $y(0) = 0$.

The potential $V (x)$ decays more quickly than $\displaystyle
\frac{l(l+1)}{x^2}$ , so for large $x$ (asymptotic region) the
Schr\"odinger equation \eqref{ode} becomes

\begin{equation}
y''(x) = \left[\frac{µl(l + 1)}{ x^{2}} + V (x)-E\right]y(x)
\end{equation}

The last equation has two linearly independent solutions
$\displaystyle k$$x$$j_{l}(kx)$ and $\displaystyle kxn_{l}(kx)$,
where $j_{l}$ and $n_{l}$ are the spherical Bessel and Neumann
functions. When $\displaystyle x \rightarrow \infty$ the solution
takes the asymptotic form

\begin{eqnarray}
\displaystyle y(x) & \approx & Akxj_{l}(kx)-Bkxn_{l}(kx) \\ &
\approx & \displaystyle D\left[sin(k x - \pi l/2) + tan(\delta_{l})
cos (kx-\pi l/2)\right]
\end{eqnarray}

where $\delta_{l}$ is called scattering phase shift and it is given
by the following expression:

\begin{equation}
\displaystyle
tan(\delta_{l})=\frac{y(x_{i})S(x_{i+1})-y(x_{i+1})S(x_{i})}{y(x_{i+1})C(x_{i})-y(x_{i})C(x_{i+1})}
\end{equation}

where $\displaystyle S(x) = k x j_{l}(k x)$, $ \displaystyle C(x) =
k x n_{l}(k x)$ and $x_{i} < x_{i+1}$ and both belong to the
asymptotic region. Given the energy we approximate the phase shift,
the accurate value of which is $\pi /2$ for the above problem. We
will use two values for the energy: $ 989.701916$ and $341.495874$.
As for the frequency $w$ we will use the suggestion of Ixaru and
Rizea in \cite{ix_ri} and \cite{ixaru85}.:

\begin{equation}
w=\left\{
\begin{array}{l}
\sqrt{E-50},\quad if \quad x\in[0,\ 6.5] \\
\sqrt{E},\qquad \quad else \quad x\in[6.5,\ 15]
\end{array}
\right.
\end{equation}

In Figure 1 we use $E = 989.701916$ and in Figure 2 we use $E =
341.495874$.

\subsubsection{Inhomogeneous Equation}

$y''=-100y+99sin(t)$, with $y(0)$=1, $y'(0)$=11, $t\in[0, 1000\pi]$.
Theoretical solution: $y(x)=\sin(t)+\sin(10t)+\cos(10t)$.\\
Estimated frequency: $w$=10.

\subsubsection{Duffing Equation}

$y''=-y-y^{3}+0.002\cos(1.01t)$, with $y(0)$=0.200426728067,
$y'(0)$=0, $t\in[0, 1000\pi]$. Theoretical solution:
$y(x)=0.200179477536\cos(1.01t)+2.46946143$
$*10^{-4}\cos(3.03t)+3.04014*10^{-7}\cos(5.05t)+3.74*10^{-10}\cos(7.07t)+...$\\
Estimated frequency: $w$=1.
\subsubsection{Nonlinear Problem}

$y''=-100y+sin(y)$, with $y(0)$=0, $y'(0)$=1, $t\in[0, 20\pi]$.The
theoretical solution is not know, but we use
$y(20\pi)=3.92823991*10^{-4}$ and $w$=10 as frequency of this
problem.

\section{Conclusions}
\label{sec2}

\begin{figure}[hbt]
\begin{center}
   \includegraphics[width=\textwidth]{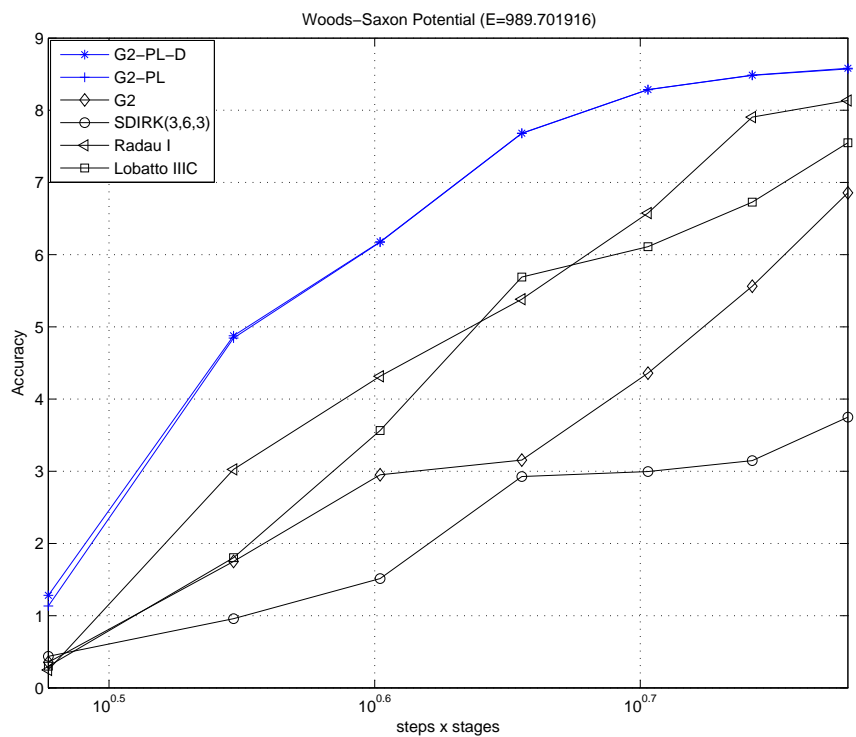}
   \caption{Resonance Problem using E=989.701916}
   \label{fig_Schrod_Eq1}
\end{center}
\end{figure}

\begin{figure}[hbt]
\begin{center}
   \includegraphics[width=\textwidth]{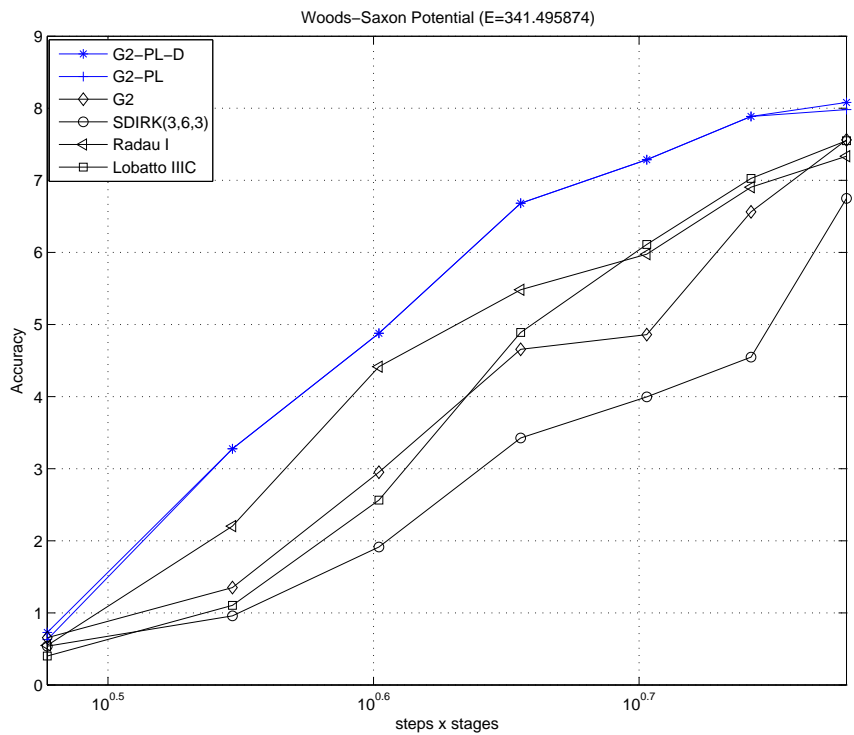}
   \caption{Resonance Problem using E=341.495874}
   \label{fig_Schrod_Eq2}
\end{center}
\end{figure}

\begin{figure}[hbt]
\begin{center}
   \includegraphics[width=\textwidth]{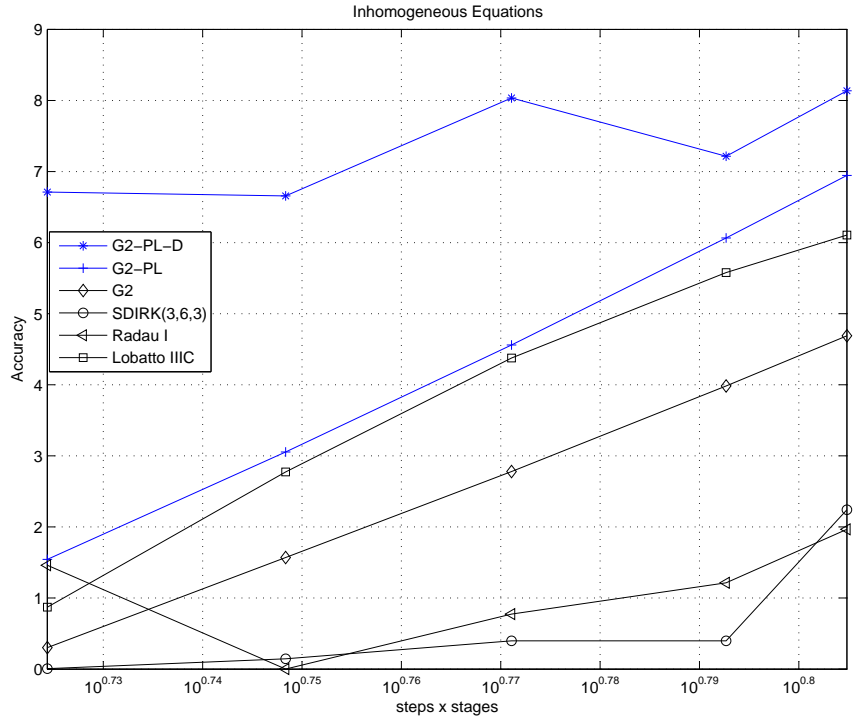}
   \caption{Inhomogeneous Equation}
   \label{fig_Inh_Eq}
\end{center}
\end{figure}

\begin{figure}[hbt]
\begin{center}
   \includegraphics[width=\textwidth]{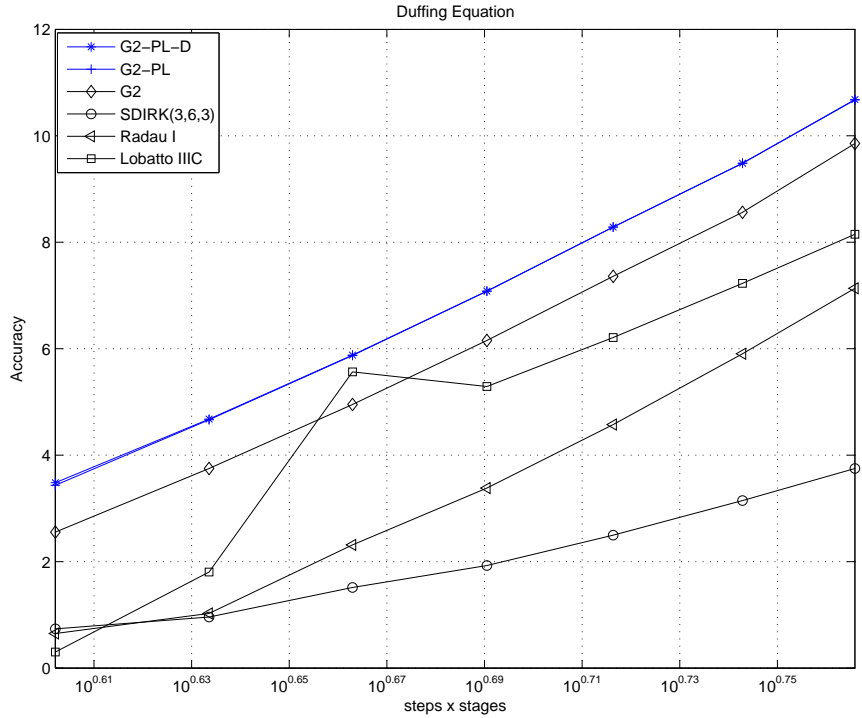}
   \caption{Duffing Equation}
   \label{fig_Duf_Eq}
\end{center}
\end{figure}

\begin{figure}[hbt]
\begin{center}
   \includegraphics[width=\textwidth]{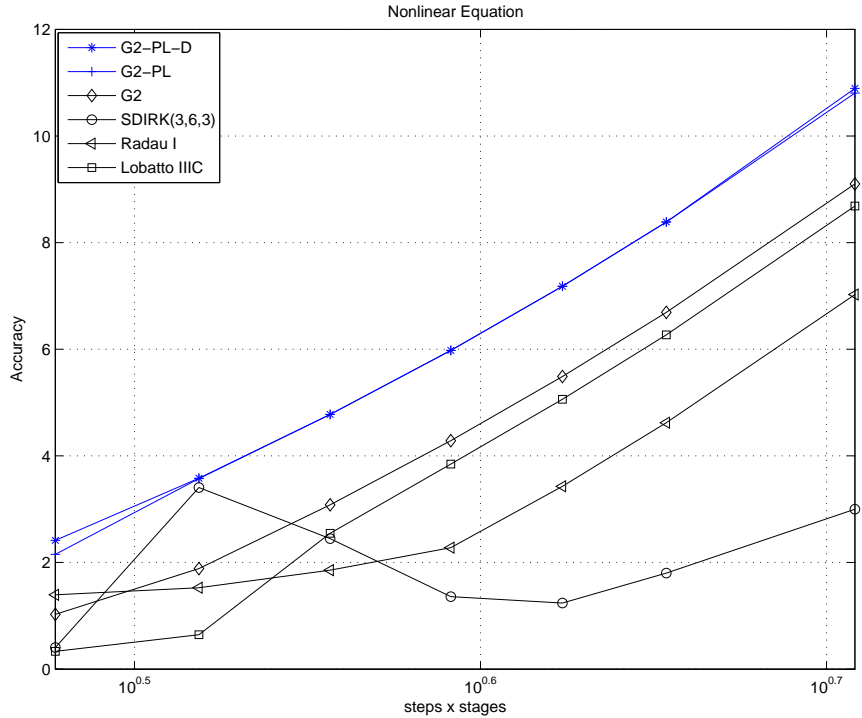}
   \caption{Nonlinear Problem}
   \label{fig_Non_Eq}
\end{center}
\end{figure}

In the following figures we present the accuracy of the tested
methods in connection with log10 of the total steps$\times$stages of
the methods. Firstly, in figure \ref{fig_Schrod_Eq1} (Resonanse
Problem) we use E=989.701916 and we observe that the second method,
which has zero phase-lag and zero dissipation, has such the same
accuracy with the first method, which has zero phased-lag, but much
better accuracy than the other methods (G2, SDIRK(3,6,3) , Radau I ,
Lobatto IIIC), while in figure \ref{fig_Schrod_Eq2} (Resonanse
Problem) we use E=341.495874 and we also see  that the second
method, which has zero phase-lag and zero dissipation, has the same
accuracy with the first method, which has zero phase-lag, but better
accuracy than the other methods (G2, SDIRK(3,6,3) , Radau I ,
Lobatto IIIC). The conclusion from the above is that the difference
in efficiency is higher when using higher energy. Secondly in figure
\ref{fig_Inh_Eq} (Inhomogeneous Equation) it can be observed that
the second method developed here with zero phase-lag and zero
dissipation has way higher efficiency than the first method
developed here with zero phase-lag, which has far higher efficiency
than the other methods G2, SDIRK(3,6,3) , Radau I , Lobatto IIIC.
Thirdly, in figure \ref{fig_Duf_Eq} (Duffing Equation) the two new
methods have almost the same efficiency but much better than the
other four methods. Finally in figure \ref{fig_Non_Eq} (Nonlinear
Equation) the two new methods have almost the same efficiency but
much better than the other four methods. This shows the importance
of zero phase-lag and zero dissipation in this type of problems.

\begin{figure}[hbt]
\begin{center}
   \includegraphics[width=\textwidth]{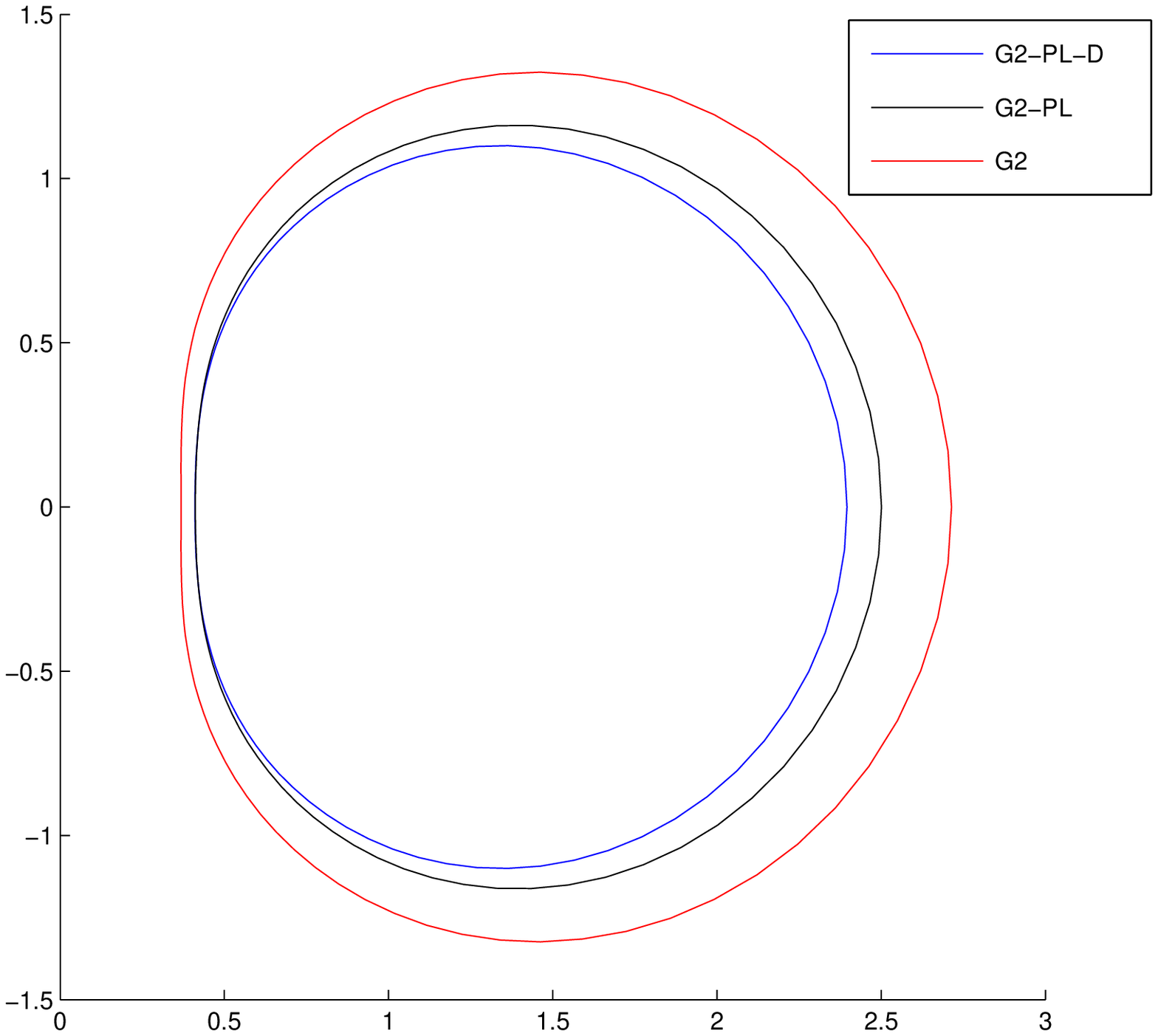}
   \caption{Stability region for the three implicit Runge-Kutta methods}
   \label{fig_stability}
\end{center}
\end{figure}

The stability region for the three methods for $\upsilon=\omega\,h=50$ are shown in Figure \ref{fig_stability}. The stability regions of the methods include the exterior of the curves. Note that the curves of all the three methods include the left half-plane.

\newpage

\textbf{Figures of the Numerical Results}

\newpage

\newpage

\newpage

\begin{thebibliography}{99}
\bibitem{rallison} D. Raptis and A.C. Allison, Exponential-fitting methods for the numerical solution of the Schr\"{o}dinger equation, Computer Physics Communications, 14, 1 (1978)

\bibitem{ix_ri} L.Gr. Ixaru, M. Rizea, A Numerov-like scheme for the numerical solution of the Schr\"{o}dinger equation in the deep continuum spectrum of energies, Comp. Phys. Comm. 19, 23-27 (1980)

\bibitem{raptis} A.D. Raptis, Exponentially-fitted solutions of the eigenvalue Schr\"{o}dinger equation with automatic error control, Computer Physics Communications, 28, 427 (1983)

\bibitem{ixaru85} L. Gr. Ixaru and M. Rizea, Comparison of some four-step methods for the numerical solution of the Schr\"odinger equation, {\it Computer Physics Communications}, {\bf 38 (3)} 329-337 (1985)

\bibitem{hou1} P. J. Van der Houwen and B. P. Sommeijer, Phase-Lag of Implicit Runge-Kutta methods, Society for Industrial and Applied Mathematics 26(1) 214-229(1989)

\bibitem{koto} Toshiyuki Koto, Phase-Lag Analysis of Diagonally Implicit Runge-Kutta Methods, Journal of Information Processing 13(3), 361-366 (1990)

\bibitem{lambert} J. D. Lambert, Numerical methods for ordinary differential systems: the initial value problem, John Wiley \& Sons, Inc., New York, NY (1991)

\bibitem{haug} E. J. Haug, D. Negrut, C. Engstler, Implicit Runge-Kutta Integration of the Equation of Multibody Dynamics in Descriptor Form, Advances in Design Automation-Proceedings of the ASME Design Autmation Conference, Sacramento, CA (1997)

\bibitem{simos2} T.E. Simos, Four-step P-stable method with minimal phase-lag, Computer Physics Communications, 115, 1 (1998)

\bibitem{hou2} P. J. Van der Houwen, B. P. Sommeijer, Diagonally Implicit Runge-Kutta methods for 3D Shallow water applications, {\em Advances in Computational Mathematics} 12, 229 - 250 (2000)

\bibitem{ref_03} T.E. Simos, P.S. Williams, A P-stable hybrid exponentially-fitted method for the numerical integration of the Schr\"odinger equation, Computer Physics Communications 131 (2000) 109-119

\bibitem{ref_28} T.E. Simos An accurate eighth order exponentially-fitted method for the efficient solution of the Schr\"odinger equation, Computer Physics Communications 125 (2000) 21-59

\bibitem{ref_27} T.E. Simos, Jesus Vigo-Aguiar, An exponentially-fitted high order method for long-term integration of periodic initial-value problems, Computer Physics Communications 140 (2001) 358-365

\bibitem{simos_aguiar_cpc} T.E. Simos, Jesus Vigo-Aguiar, A dissipative exponentially-fitted method for the numerical solution of the Schr\"odinger equation and related problems, Computer Physics Communications 152, 274-294 (2003)

\bibitem{butcher} J.C.Butcher, Numerical Methods for Ordinary Differential Equations, John Wiley \& Sons, Ltd, Chichester, England (2003)

\bibitem{hairer} E. Hairer, G. Wanner, Solving Ordinary Differential Equations II. Stiff and Differential-Algebraic Problems, Berlin Heidelberg New York: Springer-Verlag (1996)

\bibitem{anas0} Z.A. Anastassi and T.E. Simos, Special optimized Runge-Kutta methods for IVPs with oscillating solutions, {\em International Journal of Modern Physics C} {\bf 15} 1-15 (2004)

\bibitem{anastassi_newastr1} Z.A. Anastassi and T.E. Simos: A Dispersive-Fitted and Dissipative-Fitted Explicit Runge-Kutta method for the Numerical Solution of Orbital Problems, New Astronomy, 10, 31-37 (2004)

\bibitem{anastassi_newastr2} Z.A. Anastassi and T.E. Simos: A Trigonometrically-Fitted Runge-Kutta Method for the Numerical Solution of Orbital Problems, New Astronomy, 10, 301-309 (2005)

\bibitem{anastassi2} Z.A. Anastassi. and T.E. Simos, Trigonometrically Fitted Fifth Order Runge-Kutta Methods for the Numerical Solution of the Schr\"{o}dinger Equation, Mathematical and Computer Modelling, 42 (7-8), 877-886 (2005)

\bibitem{anastassi1} Z.A. Anastassi and T.E. Simos, Trigonometrically-Fitted Runge-Kutta Methods for the Numerical Solution of the Schr\"{o}dinger Equation, Journal of Mathematical Chemistry 3, 281-293 (2005)

\bibitem{anas1} Z.A. Anastassi and T.E. Simos, An optimized Runge-Kutta method for the solution of orbital problems, Journal of Computational and Applied Mathematics, 175, 1 (2005)

\bibitem{anastassi3} Z.A. Anastassi. and T.E. Simos, A Family of Exponentially-Fitted Runge-Kutta Methods with Exponential Order up to Three for the Numerical Solution of the Schr\"{o}dinger Equation, Journal of Mathematical Chemistry, 41, 1, 79-100 (2007)

\bibitem{ref_13} T.E. Simos, High-order closed Newton-Cotes trigonometrically fitted formulae for long-time integration of orbital problems, Computer Physics Communications 178, (2008), 199-207

\bibitem{franco} J. M. Franco, I. Gomez, L. Randez, SDIRK methods for stiff
ODEs with oscillating solutions, Applied Mathematics, 81, 197-209,
(1997)
\end{thebibliography}
\end{document}